\documentclass[AMA,STIX1COL]{WileyNJD-v2}

\usepackage{amsmath}
\usepackage{amssymb}

\newcommand*{\citen}[1]{%
  \begingroup
    \romannumeral-`\x 
    \setcitestyle{numbers}%
    [\cite{#1}]%
  \endgroup   
}

    \newcommand{\bfc}{\boldsymbol{c}}

    \newcommand{\bff}{\boldsymbol{f}}
    \newcommand{\bfF}{\boldsymbol{F}}

    \newcommand{\bfn}{\boldsymbol{n}}

    \newcommand{\bfx}{\boldsymbol{x}}

    \newcommand{\bftau  }{\boldsymbol{\tau  }}

\articletype{RESEARCH ARTICLE}%
\received{}
\revised{}
\accepted{}
\begin{document}

\title{Doubly Degenerate Diffuse Interface Models of Surface Diffusion}

\author[1,2]{Marco Salvalaglio*}
\author[1,2]{Axel Voigt}
\author[3]{Steven M. Wise}

\authormark{M. SALVALAGLIO \textsc{et al}}

\address[1]{\orgdiv{Institute of Scientific Computing, Department of Mathematics}, \orgname{TU Dresden}, \orgaddress{\state{01062 Dresden}, \country{Germany}}
}
\address[2]{\orgdiv{Dresden Center for Computational Materials Science}, \orgname{TU Dresden}, \orgaddress{\state{01062 Dresden}, \country{Germany}}
}
\address[3]{\orgdiv{Department of Mathematics}, \orgname{The University of Tennessee}, \orgaddress{\state{Knoxville, TN 37996}, \country{USA}}
} 

\corres{*M. Salvalaglio, TU Dresden, 01062 Dresden, Germany. \email{marco.salvalaglio@tu-dresden.de}}


    \abstract[Summary]{
We discuss two doubly degenerate Cahn-Hilliard (DDCH) models for isotropic surface diffusion. Degeneracy is introduced in both the mobility function and a restriction function associated to the chemical potential. Our computational results suggest that the restriction functions yield more accurate approximations of surface diffusion. We consider a slight generalization of a model that has appeared before, which is non-variational, meaning there is no clear energy that is dissipated along the solution trajectories. We also introduce a new variational and, more precisely, energy dissipative model, which can be related to the generalized non-variational model. For both models we use formal matched asymptotics to show the convergence to the sharp interface limit of surface diffusion.}


\maketitle

    \section{Introduction}
    
Motion by surface diffusion, where the normal velocity of a hypersurface in Euclidean space is given by the surface Laplacian of the mean curvature, plays an important role in various applications in material science. An important example is solid state dewetting. See, for example, \citen{Thompson_ARMR_2012} for a general introduction of this subject. While various direct numerical approaches have been developed for this  fourth order geometric evolution law -- see, e.g. \citen{Baenschetal_JCP_2005,Hausseretal_IFB_2005,Baoetal_JCP_2017,Barrettetal_JCP_2019} -- for many applications, diffuse interface approximations are considered, see e.g. \citen{Wise2005,Ratz2006,Yeon2006,Torabi2009,Li2009,Banas2009,Jiang2012,Salvalaglio2015a,Bergamaschini2016,Naffouti2017,Schiedung2017}. These diffuse interface approaches capture the motion of the interface implicitly as the evolution of an iso-surface of a phase field function. Typically, they are fourth-order nonlinear diffusion equations of Cahn-Hilliard type, whose solutions formally converge to those of their sharp interface counterpart, as the interface thickness tends to zero\cite{cahn1996,Ratz2006,Gugenberger2008,Dziwniketal_Non_2017}. 

The now-conventional wisdom is that, when the Cahn-Hilliard equation is paired with a polynomial free energy, for the aforementioned sharp interface convergence result to hold, it is required that the degeneracy in the so-called mobility function needs to be of sufficiently high order~\cite{Voigt2016}. As recently shown \cite{Daietal_MMS_2014,Leeetal_APL_2015,Leeetal_SIAMJAM_2016}, occasionally used second order degenerate mobility functions -- see e.g., \citen{Bhateetal_JAP_2000,Wiseetal_JCP_2006,Torabi2009,Jiang2012} -- do actually not converge to surface diffusion, but contain additional bulk diffusion contributions. This is not the case for the originally proposed phase field approximation of surface diffusion~\cite{cahn1996}, which uses a double-obstacle potential (in the deep-quench limit) instead of the double-well polynomial potential. In a series of papers \cite{dai2012,Daietal_MMS_2014,dai2016,dai2016JCP}, the authors conduct careful studies, both computational and theoretical, examining the effects on coarsening rates, asymptotic limits, and other solution properties, when using a range of diffusional mobility types in the Cahn-Hilliard equation, including one-sided and two-sided degeneracies of various orders. Though these papers did not solely target Cahn-Hilliard models for surface diffusion, they do confirm the findings in~ \citen{Daietal_MMS_2014,Leeetal_APL_2015,Leeetal_SIAMJAM_2016}.

In the diffuse interface model for surface diffusion proposed in~\citen{Ratz2006}, an additional degeneracy is introduced, following similar ideas as used for the thin film limit in classical phase field models for solidification \cite{Karma1998}. We refer to this model as a doubly degenerate Cahn-Hilliard (DDCH) equation.  This second degeneracy does not alter the asymptotic limit \cite{Ratz2006}, but actually leads to more accurate surface diffusion approximations, see, e.g., \citen{Backofen2019}. In fact several simulations for realistic applications in materials science consider this additional degeneracy \cite{Albani2016,Salvalaglio2015,SalvalaglioPRB2016,Salvalaglio2017a,Salvalaglio2017b,Naffouti2017,Geslin2019,Albani2019}. However, the model of \citen{Ratz2006} also has a drawback. It is non-variational, that is, there is no known free energy that is dissipated along solution trajectories. This makes it harder to prove properties of solutions and derive certain numerical stabilities. Furthermore, it excludes variational derivations in more complex multi-physics applications involving surface diffusion.  

This paper is organized as follows. In Section~\ref{sec:2}, we reintroduce the DDCH model of \citen{Ratz2006}, generalizing it in a trivial way and recalling some asymptotic convergence results. In Section~\ref{sec:3}, we introduce a new DDCH variational model for surface diffusion, and we analyze some of its properties. We make connections in Section~\ref{sec:3} between the new DDCH variational model and (i) the non-variational DDCH model of \citen{Ratz2006} and (ii) a deGennes diffuse interface model from another modeling context. We provide a simple, but powerful, numerical integration scheme in Section \ref{sec:4}, and we compare the numerical results of the  models in Section \ref{sec:5} with various parameter choices. The matched asymptotic analysis for the (new) variational and the (well-used) non-variational surface diffusion models are provided in Appendices~\ref{sec-asymptotics-model-v} and \ref{sec-asymptotics-model-nv}, respectively.
    
    \section{A Non-Variational Diffuse Interface Model} 
    \label{sec:2}
    
In this section we reintroduce a non-variational DDCH model that approximates sharp interface motion by surface diffusion. Suppose that $\Omega$ is a bounded subset of $\mathbb{R}^d$. Let $u:\Omega \to \mathbb{R}$ be an order/phase field parameter. Consider the following Cahn-Hilliard-type model proposed in \citen{Ratz2006}:
    \begin{align}
\partial_t u & = \frac{1}{\varepsilon} \nabla\cdot\left(M_0( u)\nabla w \right),
    \label{eqn-ch-G-1}
    \\
G_0( u) w & = \frac{1}{\varepsilon}f'( u) - \varepsilon  \Delta  u,
    \label{eqn-ch-G-2}
    \end{align}
where $\varepsilon >0$ is a small parameter (relative to the domain size) that is related to the thickness of the diffuse interface. $M_0$ denotes the mobility function, which is defined as
    \[
M_0(u) = \mu u^2(1-u)^2, \quad \mu = 36,
    \]
and $f$ is the quartic, symmetric double well potential function
    \[
f(u) = \frac{\omega}{4}u^2(1-u)^2, \quad \omega = 72.    
    \]
Its minima, 0 and 1, represent the pure phase states. Let us assume that, for simplicity, $ u$, the phase field, and $w$, called the chemical potential, satisfy periodic boundary conditions with respect to $\Omega$. Clearly, mass is conserved in the system: $d_t \int_\Omega  u(\bfx,t)\, d\bfx = 0$. 

To keep the model as general as possible, let us assume that $G_0$ -- which we shall call the \emph{diffusion restriction function} -- is defined as
    \[
G_0(u) := \eta |u|^p|1-u|^p, \quad p \ge 0, \quad \eta > 0.
    \]
The authors of \citen{Ratz2006} assumed that $p=2$ and $\eta = 30$, and they showed using matched asymptotic analysis that, as $\varepsilon\searrow 0$, solutions converge formally to those of a sharp interface model of surface diffusion. Let $\Sigma$ define a hypersurface that coincides with the 0.5 level set of $ u$ at some time $t$. We say that $\Sigma$ evolves by surface diffusion if 
    \begin{equation}
v = \Delta_\Sigma \kappa,
    \label{eqn-surf-diff}
    \end{equation}
where $v$ is the normal velocity at the sharp interface $\Sigma$, $\Delta_\Sigma$ is the surface Laplacian, and $\kappa$ is the mean curvature of $\Sigma$ \cite{Mullins1957,Mullins1959}. We show in Appendix~\ref{sec-asymptotics-model-nv}, by a simple calculation, that for the generalized model, upon taking 
    \begin{equation}
\eta = \eta_\star(p) :=  \frac{\Gamma(2+2p)}{\left(\Gamma(1+p)\right)^2}, \quad p \ge 0,
    \label{eqn-eta-normalization}
    \end{equation}
where $\Gamma$ is the usual (Bernoulli) Gamma function, the limiting law as $\varepsilon\searrow 0$ is again motion by surface diffusion \eqref{eqn-surf-diff}. The normalization in \eqref{eqn-eta-normalization} is called the \emph{diffusion normalization}. (Note that $\eta_\star(p=2) = 30$, which recovers the result in~\citen{Ratz2006}.)

The model \eqref{eqn-ch-G-1} -- \eqref{eqn-ch-G-2} is a type of doubly degenerate Cahn-Hilliard (DDCH) equation, because degeneracies are associated with both the mobility and the restriction functions.   We refer to this model \eqref{eqn-ch-G-1} -- \eqref{eqn-ch-G-2} as the \textit{non-variational DDCH model} (or, Model NV, for short), because it is not arrived at through a variational principle. Indeed, it is not clear whether there is a related energy that is dissipated along the solution trajectories of \eqref{eqn-ch-G-1} -- \eqref{eqn-ch-G-2}. In other words, it seems that the model is not free energy dissipative, or thermodynamically consistent.  Model NV becomes degenerate as $u$ approaches the pure state values $ u=0$ and $ u=1$. Indeed, both the mobility $M_0$ and the restriction function $G_0$ are degenerate when $ u = 0,1$. There is ample numerical evidence to suggest that there is an important benefit owing to these degeneracies, namely, solutions to Model NV remain in the physically relevant region $0<  u < 1$, as long as the initial data have this feature. This, to our knowledge, has yet to be verified theoretically. 

Incidentally, this property -- $0 <  u(\, \cdot \, , 0) < 1 \implies 0 <  u(\, \cdot \, , t) < 1$, for all $t\ge 0$ -- is often referred to as a \emph{positivity preserving} property. To see why, note that for a typical binary model, like the traditional Cahn-Hilliard equation, the concentrations of the species $A$ and $B$ are given by $ u_A =  u$ and $ u_B = 1- u$. Therefore we have the positivity of the concentrations, $ u_A >0$ and $ u_B >0$, iff $0 <  u < 1$. We use this terminology herein.

We point out that one can regularize Model NV so that it is defined for all values of $ u$, such that the asymptotic limit should remain unchanged. This may be important for numerical implementation, since most numerical schemes cannot guarantee that solutions remain positive, even when this property is known to hold for the PDE. For example, one can reintroduce the mobility as
	\[
M_\alpha(u) = \mu u^2(1-u)^2 + \alpha\varepsilon, \quad \mu = 36, \quad \alpha \ge 0.
    \]
To make the restriction function defined, continuous, positive, and differentiable on all of $\mathbb{R}$, one can regularize as follows:
    \[
G_\alpha(u) = \sqrt{\eta^2 \left(u^2(1-u)^2 \right)^p + \alpha^2\varepsilon^2} , \quad p \ge 0, \quad \alpha \ge 0.  
    \]
Thus, setting $\alpha = 0$ we obtain the model above.

Upon choosing $p = 0$, $\eta = 1$, and  $\alpha = 0$, one obtains a more familiar (singly) degenerate Cahn-Hilliard (DCH) equation:
    \begin{align}
\partial_t u & = \frac{1}{\varepsilon} \nabla\cdot\left(M_0( u)\nabla w \right),
    \label{eqn-ch-G-1-b}
    \\
w & = \frac{1}{\varepsilon}f'( u) - \varepsilon  \Delta  u.
    \label{eqn-ch-G-2-b}
    \end{align}
Matched asymptotic analysis shows that solutions of this model also formally converges to those of a sharp interface model of surface diffusion, as $\varepsilon\searrow 0$. However, quite importantly, numerical results suggest such diffuse interface solutions converge more slowly to the sharp interface model compared with those of Model NV, \eqref{eqn-ch-G-1} -- \eqref{eqn-ch-G-2}. See, for example, \citen{Backofen2019}. The computational results in Section~\ref{sec:5} support this point further. However, model \eqref{eqn-ch-G-1-b} -- \eqref{eqn-ch-G-2-b} is free energy dissipative, that is, thermodynamically consistent. Formally, solutions to this model dissipate the free energy
    \[
F[ u] = \int_\Omega \left\{ \frac{1}{\varepsilon} f( u)+\frac{\varepsilon}{2}|\nabla u|^2 \right\} d\bfx 
    \]
at the rate
    \[
d_t F[ u] = - \frac{1}{\varepsilon}\int_{\Omega} M_0( u) |\nabla w|^2\, d\bfx .
    \]

    \section{A Variational Diffuse Interface Model} 
    \label{sec:3}
    
Is there an energy dissipative DDCH model that is related to Model NV? The answer is a qualified \emph{yes}. In this section, we show that there is a variational model such that Model NV is its approximation, in certain cases. To this end, consider the free energy 
	\begin{equation}
F[u] = \int_\Omega g_0(u) \left(\frac{1}{\varepsilon}f(u) +\frac{\varepsilon}{2}|\nabla u|^2 \right)d{\bf x} ,
    \label{eqn-reg-energy}
    \end{equation}
where $f$ is the same quartic potential as before. Here we assume that $g_0$ is a singular function of the form
    \[
g_0(u) = \frac{1}{\gamma |u|^p |1-u|^p}, \quad p \ge 0, \quad \gamma >0.
    \]
We call $g_0$ the \emph{energy restriction} function. If necessary, $g_0$ can be regularized so that it is defined, continuous, and differentiable for all $u$: 
	\[
g_\alpha(u) = \frac{1}{\sqrt{\gamma^2{(u^2(1-u)^2)}^p + \alpha^2\varepsilon^2}}, \quad p\ge 0, \quad \alpha \ge 0.
    \]
Thus, setting $\alpha = 0$, we obtain the function above.

We are interested in the energy dissipative flow
	\begin{align}
\partial_t u &= \frac{1}{\varepsilon}\nabla\cdot \left(M_0 (u)\nabla w \right),
    \label{eqn-ch-little-g-1}
    \\
w &= g_0(u)\frac{f'(u)}{\varepsilon} -\varepsilon \nabla\cdot\left(g_0(u)\nabla u \right) + g_0'(u)\left(\frac{1}{\varepsilon}f(u) +\frac{\varepsilon}{2}|\nabla u|^2 \right) ,
    \label{eqn-ch-little-g-2}
	\end{align}
where $w = \delta_u F$ is the chemical potential and $\delta_u F$ is the variational derivative with respect to $u$. Here, for simplicity, we have assumed natural or periodic boundary conditions on $\Omega$.  System \eqref{eqn-ch-little-g-1} -- \eqref{eqn-ch-little-g-2} is another type of doubly degenerate Cahn-Hilliard (DDCH) equation, which we refer to as the \emph{variational DDCH model} (or, Model V, for short). 

The derivative $g'_0$ is, of course, singular,
    \[
 g'_0(u)=p\frac{2u  - 1}{\gamma u^{p+1}(1-u)^{p+1} } , \quad 0 < u<1 ,
    \]
which, it seems, help to keep solutions positive. The regularized version is non-singular,
	\[
g'_\alpha(u) =  \frac{ p\left[u^2(1-u)^2\right]^{p-1}u(1-u)(2u-1)}{\gamma \left(\left[(u^2(1-u)^2)\right]^p + \alpha^2 \varepsilon^2\right)^{3/2} } ,
    \]
though its values at 0 and 1 become increasing large as $\alpha\searrow 0$. There are some open analysis questions related to the validity of the model when the regularization vanishes ($\alpha=0$). For example, do solutions have the positivity preserving property? Do weak solutions exist and are they regular, in some sense? For what values of $p$ and $\gamma$ does the model make sense? Our early simulation results -- as well as some early theoretical results, not reported here -- seem to support the validity of a positivity preserving property, when $0 < p < 2$.  In any case, formally, it is clear that the energy $F$ given in \eqref{eqn-reg-energy} is dissipated along the solution trajectories of the dynamical system \eqref{eqn-ch-little-g-1} -- \eqref{eqn-ch-little-g-2}, that is, system \eqref{eqn-ch-little-g-1} -- \eqref{eqn-ch-little-g-2} is a free energy dissipative dynamical system.

    \subsection{Relation to Model NV}

Let us  now see how the variational/dissipative DDCH model (Model V) may be related to the non-variational DDCH model (Model NV).  Assume that, to leading order in $\varepsilon$, the profile of the phase field function $u$ perpendicular to the  diffuse interface is the hyperbolic tangent function
	\begin{equation}
u(\bfx,t) \approx u_I(z) = \frac{1}{2}\left(1+\tanh\left(\frac{3 z}{\varepsilon} \right)\right), 
    \label{eq:tanh}
    \end{equation}
where $z$ is the coordinate perpendicular to the interface, a type of  interface distance function. Then, we have the asymptotic approximation, to leading order in $\varepsilon$,
	\begin{equation}
\frac{1}{\varepsilon}f(u) \approx \frac{\varepsilon}{2}|\nabla u|^2 \quad \mbox{or} \quad \frac{1}{\varepsilon}f(u) - \frac{\varepsilon}{2}|\nabla u|^2 = \mathcal{O}(\varepsilon).
    \label{approx-tanh}
    \end{equation}
(We will justify these approximations in the appendices, though they are somewhat standard.) Since 
    \[
\nabla \cdot \left[ g_0(u) \nabla u \right]=g_0\Delta u + g_0'(u)|\nabla u|^2,
    \]
we can approximate the chemical potential as
	\[
w \approx g_0(u) \frac{1}{\varepsilon}f'(u) - \varepsilon g_0(u) \Delta u .
    \]
In other words,
    \[
G_0(u) w = \gamma |u|^p|1-u|^p w \approx \frac{1}{\varepsilon}f'(u) - \varepsilon  \Delta u , \quad \eta = \gamma.
    \]
Thus, Model V is related to Model NV through approximation, if the interfacial profile is a hyperbolic tangent, and if $\eta =\gamma$. In the next section, we find that there is a logical choice for $\gamma$, given $p$, for Model V, the so-called energy normalization.
    
    \subsection{Finite Energy and Energy Normalization}    
    
When, if ever, is the energy \eqref{eqn-reg-energy} of a diffuse interface solution with the hyperbolic tangent profile finite? To answer this, let us work in two space dimensions, for simplicity.  If \eqref{eq:tanh} is valid to leading order, then, using the approximation~\eqref{approx-tanh} the energy can be approximated as 
	\[
F[u] \approx \frac{2}{\varepsilon}\int_\Omega g_0(u) f(u) \, d{\bf x} \approx |\Sigma|\frac{2}{\varepsilon} \int_{-\infty}^\infty g_0(u_I(z)) f(u_I(z)) \, d z ,
    \]
where $|\Sigma|$ is the length of the diffuse interface, measured in the tangential direction along the $u = 0.5$ level curve. (In three dimensions, the surface area of the $0.5$ iso-surface would appear.) For $p=1$, taking $\gamma = 6$, we find that the energy is finite and
	\[
F[u] \approx |\Sigma|\frac{2}{\varepsilon} \int_{-\infty}^\infty \frac{18 u_I^2(1-u_I)^2}{6u_I(1-u_I)} dz = |\Sigma|\frac{6}{\varepsilon} \int_{-\infty}^\infty u_I(1-u_I) \, dz = |\Sigma|.
    \]
This is the usual measure of energy for isotropic motion by mean curvature in two dimensions, namely, the energy is proportional to the interface length.  In fact, we can generalize the last result: for any $p\in[0,2)$,
	\begin{align*}
F[u] &\approx |\Sigma|\frac{2}{\varepsilon} \int_{-\infty}^\infty \frac{18 u_I^2(1-u_I)^2}{\gamma |u_I|^p|1-u_I|^p}\, dz
    \\
&= |\Sigma|\frac{36}{\varepsilon\gamma } \int_{-\infty}^\infty |u_I|^{2-p}|1-u_I|^{2-p}\,  dz 
    \\
&= \frac{6|\Sigma|}{\gamma} \int_{u_I = 0}^{u_I=1}|u_I|^{1-p}|1-u_I|^{1-p}\, du_I \qquad (\mbox{using} \ du_I = \frac{6}{\varepsilon}u_I(1-u_I)\, dz)
    \\
&= \frac{6|\Sigma|}{\gamma} \frac{(\Gamma(2-p))^2}{\Gamma(4-2p)} .
    \end{align*}
We can use this calculation to pick the value of $\gamma$ in the general case so that $F[u]\approx |\Sigma|$, for all $p\in [0,2)$. For Model V, let us choose
    \begin{equation}
\gamma = \gamma_\star(p) := 6\frac{(\Gamma(2-p))^2}{\Gamma(4-2p)}, \quad p \in [0,2).
    \label{eqn-energy-normalize-gamma}
    \end{equation}
We refer to this choice as the \emph{energy normalization}. This is different from the diffusion normalization~\eqref{eqn-eta-normalization}, though the two values are equal when $p=0$ and $p=1$:
    \[
 1 = \gamma_\star(0) = \eta_\star(0)  \quad \mbox{and} \quad 6 = \gamma_\star(1) = \eta_\star(1). 
    \]
In any case, observe that,   
    \[
\lim_{p \nearrow 2}\gamma_\star(p) = +\infty.
    \]
Consequently, the energy cannot be made finite for the hyperbolic tangent profile when $p \ge 2$. This seems to imply that  Model V cannot be made sensible for $p \ge 2$.

For $0\le p <2$, using the energy normalization, we show in Appendix~\ref{sec-asymptotics-model-v} that solutions to Model V, \eqref{eqn-ch-little-g-1} -- \eqref{eqn-ch-little-g-2},  converge, as $\varepsilon \searrow 0$, to the solutions of the sharp interface model of surface diffusion~\eqref{eqn-surf-diff}.  This fact leads us to the following conclusion: for $p=1$, choosing the energy normalization in Model V and the diffusion normalization in Model NV, the two models are directly connected via the approximation outlined in the previous section. We will demonstrate this connection in the computational results of Section~\ref{sec:4}.  For now, the case $p=1$ has other interesting properties that we will explore in the next section.

    \subsection{Relation to a deGennes Model}

Suppose in Model V that we select the parameters
	\[
\omega = 72, \quad p = 1, \quad  \alpha = 0, \quad \mbox{and} \quad \gamma = 6.	
	\]
Under the assumption that solutions are restricted to the physically meaningful region of phase space $0 \le  u  \le  1$, the energy \eqref{eqn-reg-energy} is formally equivalent to 
	\begin{equation}
{\mathcal F}_0[u] = \int_\Omega \left\{ I(u) +  \frac{3}{\varepsilon}u(1-u) +\frac{\varepsilon}{2}g_0(u)|\nabla u |^2  \right\}d\bfx ,
    \label{eqn-energy-indicator}
	\end{equation}
where $I$ is the obstacle potential
	\[
I(u) = 
	\begin{cases}
0 & \mbox{if} \quad  0\le u \le 1
	\\
+\infty  & \mbox{if} \quad u < 0 \quad \mbox{or} \quad u>1
	\end{cases}	
.
	\]
The energy density $\frac{\varepsilon}{2}g_0(u)|\nabla u |^2$ is called the \emph{deGennes gradient energy density}~\cite{dong2019,li2017}. Here, it seems, the obstacle potential does not play an essential role, since the singular deGennes gradient energy density is expected to keep the solution within the physically meaningful range in a gradient flow setting.  However, this version of the energy, ${\mathcal F}_0$, is instructive, as it can be viewed as the deep-quench limit of the following Flory-Huggins-deGennes-type model:
	\begin{equation}
{\mathcal F}_\theta[u] = \int_\Omega \left\{ \theta \left( u\log(u) + (1-u) \log(1-u)\right)  +  \frac{3}{\varepsilon}u(1-u) +\frac{\varepsilon}{2}g_0(u)|\nabla u |^2  \right\}d\bfx .
    \label{eqn-energy-Florry}
	\end{equation}
In other words, ${\mathcal F}_0$ is obtained in the (zero temperature) limit $\theta \searrow 0$ of ${\mathcal F}_\theta$. We note that ${\mathcal F}_\theta$ is a model that is commonly used in the study of polymer or hydrogel materials~\cite{dong2019,li2017}. Furthermore, ${\mathcal F}_\theta$ is also reminiscent of the standard Flory-Huggins free energy considered in~\citen{cahn1996}, before the deep-quench assumption is invoked.

One interesting feature of the deep-quench energy, ${\mathcal F}_0$, is that hyperbolic tangent solutions are one-dimensional minimizers. This implies that, for the associated conserved and dissipative dynamical system, evolving diffuse interfaces will have the hyperbolic tangent shape. This is in contrast to the deep-quench model considered in~\citen{cahn1996}, where diffuse interfaces have cosine-like profiles. 

In future works, we will address some theoretical issues related to Model V, for example, the existence of weak solutions~\cite{dai2016}, the validity of the positivity-preserving properties of the PDE and associated numerical schemes, et cetera. For the moment, we will explore some numerical solutions of Models V and NV, comparing their characteristics and properties, under the assumptions that the PDE solutions are sufficiently well-behaved.

    \section{Integration Schemes}
    \label{sec:4}

To approximate solutions of Models NV and V, we employ the AMDiS finite element software package~\cite{Vey2007,Witkowski2015}, which allows for adaptive mesh refinement around the evolving diffuse interface. To remove some of the stiffness, a linear implicit-explicit (IMEX) integration scheme is used to approximate solutions of Eqs.~\eqref{eqn-ch-little-g-1} -- \eqref{eqn-ch-little-g-2}:
    \begin{equation}
    \begin{split}
\frac{u^{n+1}}{\tau_n} - \frac{1}{\varepsilon}\nabla \cdot \left[ M_\alpha(u^{n}) \nabla w^{n+1}\right]&= \frac{u^{n}}{\tau_n} ,
    \\
w^{n+1}+\varepsilon\nabla \cdot \left[ g_\alpha(u^{n}) \nabla u^{n+1}\right]-\frac{1}{\varepsilon}\left[r(u^{n})+\chi s(u^{n})\right]u^{n+1}+& 
\\
+\varepsilon g_\alpha'(u^{n})\left( \frac{\chi}{2}-1\right) \nabla u^{n} \cdot \nabla u^{n+1}&=
q(u_n),
    \\
    \end{split}
    \label{eq:system_full}
    \end{equation}
where
    \begin{align*}
q(u_n)&=\frac{1}{\varepsilon}g_\alpha(u^{n})f'(u^{n})+\frac{\chi}{\varepsilon}g'_\alpha(u^{n})f(u^{n})-\frac{1}{\varepsilon}\left[r(u^{n})+\chi s(u^{n})\right]u^{n}
    \\
r(u^{n})&=[g'_\alpha(u^{n})f'(u^{n})+g_\alpha(u^{n})f''(u^{n})]
    \\
s(u^{n})&=[g'_\alpha(u^{n})f'(u^{n})+g''_\alpha(u^{n})f(u^{n})] 
    \end{align*}
account for the linearizations of the $g_\alpha(u^{n+1})f'(u^{n+1})$ and $g'_\alpha(u^{n+1})f(u^{n+1})$ terms around $u^{n}$. The integer $n$ is the time step index; $\tau_n>0$ is the time stepsize at step $n$; and $u^0$ is the initial condition. We use a very small regularization constant, $\alpha>0$, as the current numerical scheme is not designed to preserve the expected positivity of the solution. Positivity-preserving schemes are being developed and will be reported on in later works. $\chi$ is an auxiliary variable, such that $\chi=1$ selects Model V, while $\chi=0$ enforces the approximation $f(u)/\varepsilon \approx (\varepsilon/2)|\nabla u|^2$ and allows for the integration of Model NV. Indeed, with $\chi=0$ the system can be rewritten as
    \begin{equation}
    \begin{split}
\frac{u^{n+1}}{\tau_n} - \frac{1}{\varepsilon}\nabla \cdot \left[ M_\varepsilon(u^{n}) \nabla w^{n+1}\right]=& \frac{u^{n}}{\tau_n},
    \\
G_\alpha(u^{n})w^{n+1}+\varepsilon \Delta u^{n+1}-\frac{1}{\varepsilon}f''(u^{n})u^{n+1}=& \frac{1}{\varepsilon}f'(u^{n})-\frac{1}{\varepsilon}f''(u^{n})u^{n},
    \end{split}
    \label{eq:system_app}
    \end{equation}
where $G_\alpha(u^{n})=1/g_\alpha(u^{n})$, with $p=1$ and $\gamma = \eta = 6$. We thus obtain a semi-implicit integration scheme for Eqs.~\eqref{eqn-ch-G-1-b} -- \eqref{eqn-ch-G-2-b}. In the following, we use the system \eqref{eq:system_full} with $\chi=1$ to integrate Model V and the system \eqref{eq:system_app} to integrate Model NV. Notice that this distinction strictly holds true for $p>0$. Adaptive mesh refinement is exploited with a fine spatial discretization $h_{\rm in}$ within the diffuse interface, namely where $0.05<u<0.095$, and a coarse spatial discretization $h_{\rm out}$ elsewhere. The former is scaled with $\varepsilon$ to ensure the same number of elements within the diffuse interface for any $\varepsilon$, while the latter is fixed. The model and numerical parameters are set as follows: $\alpha=10^{-4}$, $\tau_n=\tau=10^{-5}\varepsilon$, $h_{\rm in}=\varepsilon/10$, $h_{\rm out}=0.125$.

The results of the numerical integration of the DDCH equations are compared with the dynamics of the corresponding sharp-interface (SI) evolution in two space dimensions. In particular, we will consider the SI evolution of the closed curve, $\Sigma$, corresponding to the $0.5$ level-set of the initial condition data for the DDCH models. This curve is discretized by a uniform set of points $\{ \mathbf{x}_i \}$. The discrete evolution law in terms of the velocities directed along the outward-pointing unit normals $\hat{\mathbf{n}}$ read 
    \begin{equation}
    \begin{split}
\frac{\partial \mathbf{x}_i}{\partial t} = \hat{\mathbf{n}}_i \tilde\Delta_{\Sigma} \kappa_i ,
    \end{split}
    \label{eq:sharpinterface}
    \end{equation}
where $\tilde\Delta_{\Sigma}$ is a finite difference approximation of the surface Laplacian (in this reduced dimension setting) and $\kappa_i$ is the local curvature, computed as $1/r_i$ with $r_i$ the radius of curvature. The latter is approximated as the radius of the almost osculating circle, i.e., the circle passing through $\mathbf{x}_{i}$ and $\mathbf{x}_{i \pm 1}$. We note that several other methods have been proposed and used instead of the approximation in Eq.~\eqref{eq:sharpinterface} (see, e.g, \citen{Baenschetal_JCP_2005,Hausseretal_IFB_2005,HausserJSC2007,BarrettNUMMAT2008,WangPRB2015,Baoetal_JCP_2017}). Here, according to the simple evolution it has been sufficient to use a simple finite-difference scheme exploiting a forward Euler discretization of time, with a sufficiently small time stepsize.

    \section{Numerical Results} \label{sec:5}

We numerically compare Model V, Model NV and the standard DCH model ($p=0$) with the sharp-interface (SI) solution for surface diffusion. 

    \begin{figure}[h]
    \center
\includegraphics*[width = 0.75 \textwidth ]{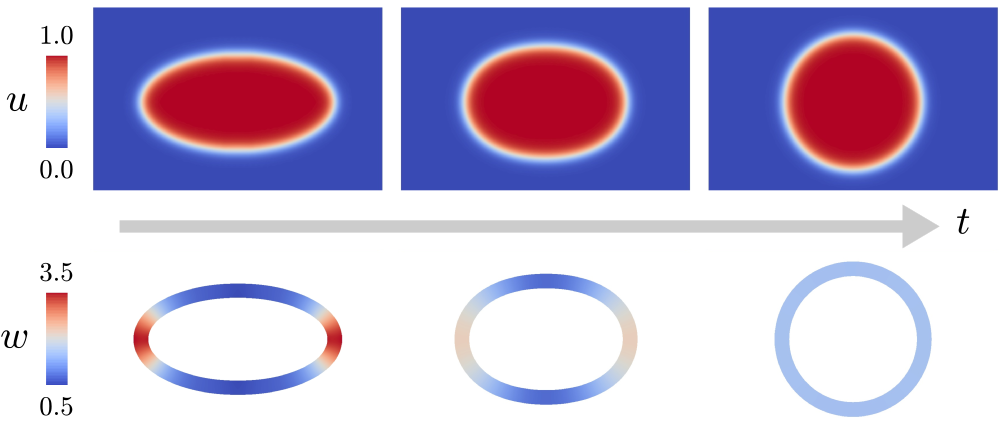} 
\caption{Evolution obtained by Model V (Eqs.~\eqref{eqn-ch-little-g-1} -- \eqref{eqn-ch-little-g-2} discretized as in \eqref{eq:system_full}) of an ellipse with semi-axis along x axis ($A_x = 1.0$) and along the y axis ($A_y = 0.5$). Illustrated by means of $u$ (top) and $w$ (bottom), the latter shown in the region where $0.1<u<0.9$. $\varepsilon=0.2$, $p=1$.}
    \label{fig:figure1}
    \end{figure}

\begin{figure}[h]
    \center
\includegraphics*[width = \textwidth ]{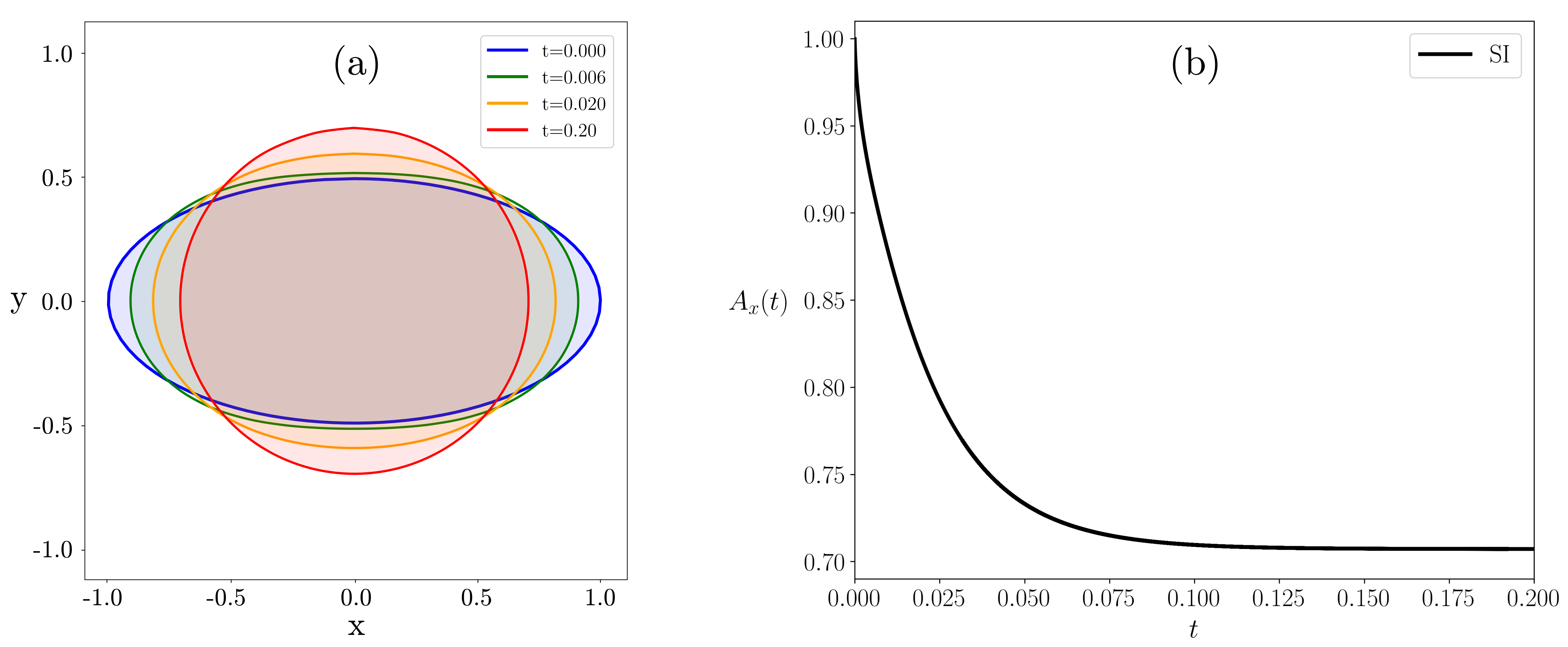} 
\caption{Sharp interface evolution (Eq.~\eqref{eq:sharpinterface}) of an ellipse with semi-axis as in Figure~\ref{fig:figure1}. (a) Representative steps during the evolution. (b) $A_x(t)$.}
    \label{fig:figure2}
    \end{figure}

Figure~\ref{fig:figure1} illustrates the evolution of an ellipse obtained by integrating Model V using \eqref{eq:system_full}. The initial condition $u^0$ is set by evaluating Eq.~\eqref{eq:tanh} with $r$ the signed distance to the considered ellipse. Figure~\ref{fig:figure2} illustrates the evolution of the $u=0.5$ level-set by the SI approach, i.e. by Eq.~\eqref{eq:sharpinterface}. A few representative profiles during the evolution are shown in Figure~\ref{fig:figure2}(a), while the change over time of the x semi-axis $A_x(t)$ is shown in Figure~\ref{fig:figure2}(b).

           \begin{figure}[h]
    \center
\includegraphics*[width = \textwidth ]{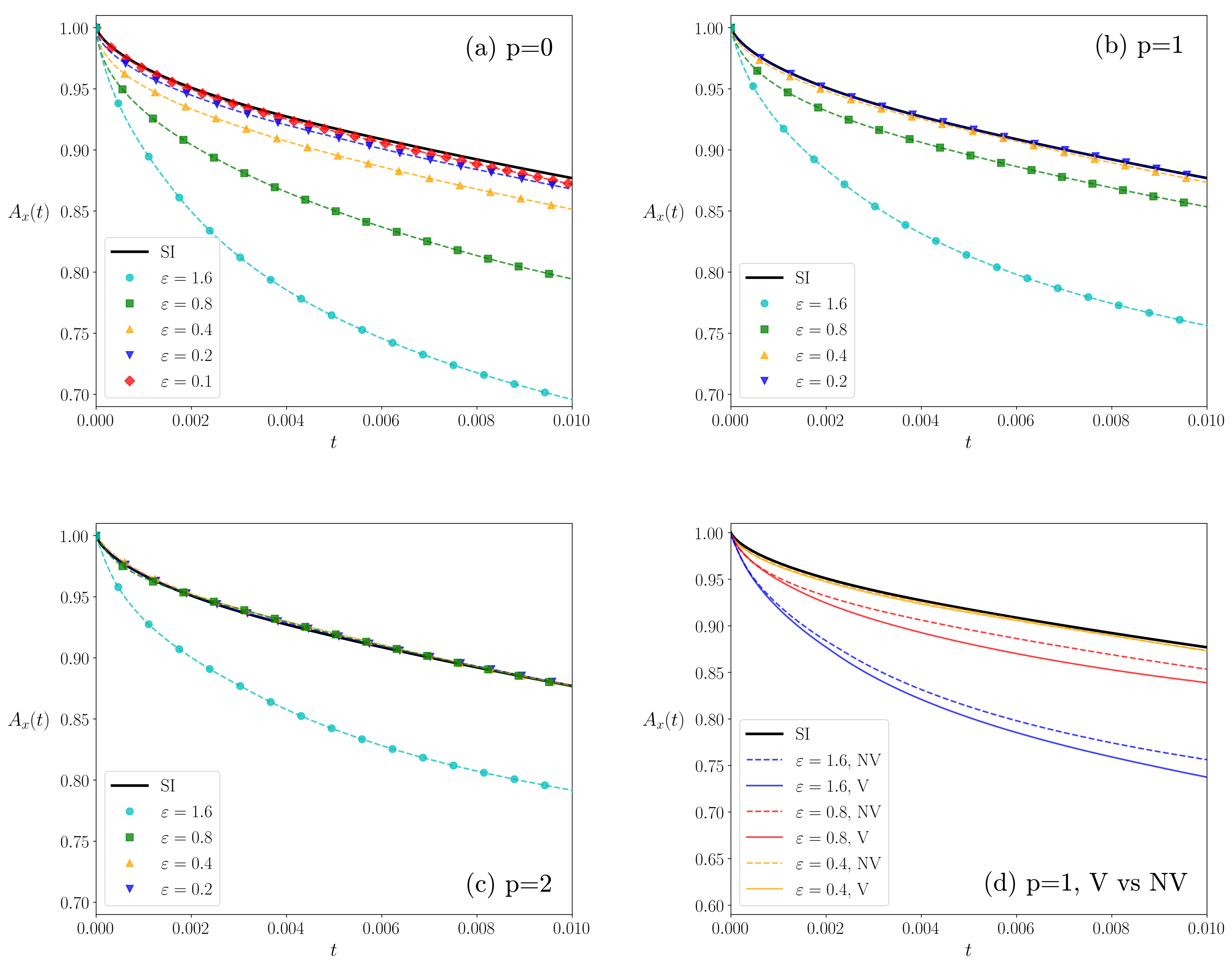} 
\caption{Results of the integration of DDCH models illustrated by means of $A_x(t)$ of the ellipse in Figure~\ref{fig:figure1} for different values of $\varepsilon$ and $p$. (a)-(c) Model NV with $p=0$, $p=1$ and $p=2$, respectively. Different curves correspond to different values for $\varepsilon$. (d) Comparison between Model V and Model NV for $p=1$. The evolution obtained by SI is reported as reference in all panels.}
    \label{fig:figureP}
    \end{figure}
    
Figure \ref{fig:figureP} shows the behavior of the DDCH models for different values of $p$ and different $\varepsilon$. We recall that a linear scaling of the timestep and of the mesh size with $\varepsilon$ have been considered. The approximation properties of the DDCH models to the sharp interface limit of surface diffusion are summarized. While all DDCH models converge to the SI solution, the error is lower for  $p=1$ and $p=2$ compared to the standard DCH model ($p=0$). For $p=1$ both Model V and Model NV can be considered. Convergence to the SI limit is similar in these cases, but we obtain slightly lower errors with Model NV (see Figure~\ref{fig:figureP}(d)). This may be ascribed to the larger number of operators to be considered for the integration of Model V. We recall that for $p=1$ Model V coincides with Model NV under the assumption $(1/\varepsilon)f(u)\approx (\varepsilon/2)|\nabla u|^2$.

\begin{figure}[h]
\center
\includegraphics*[width = 0.7 \textwidth ]{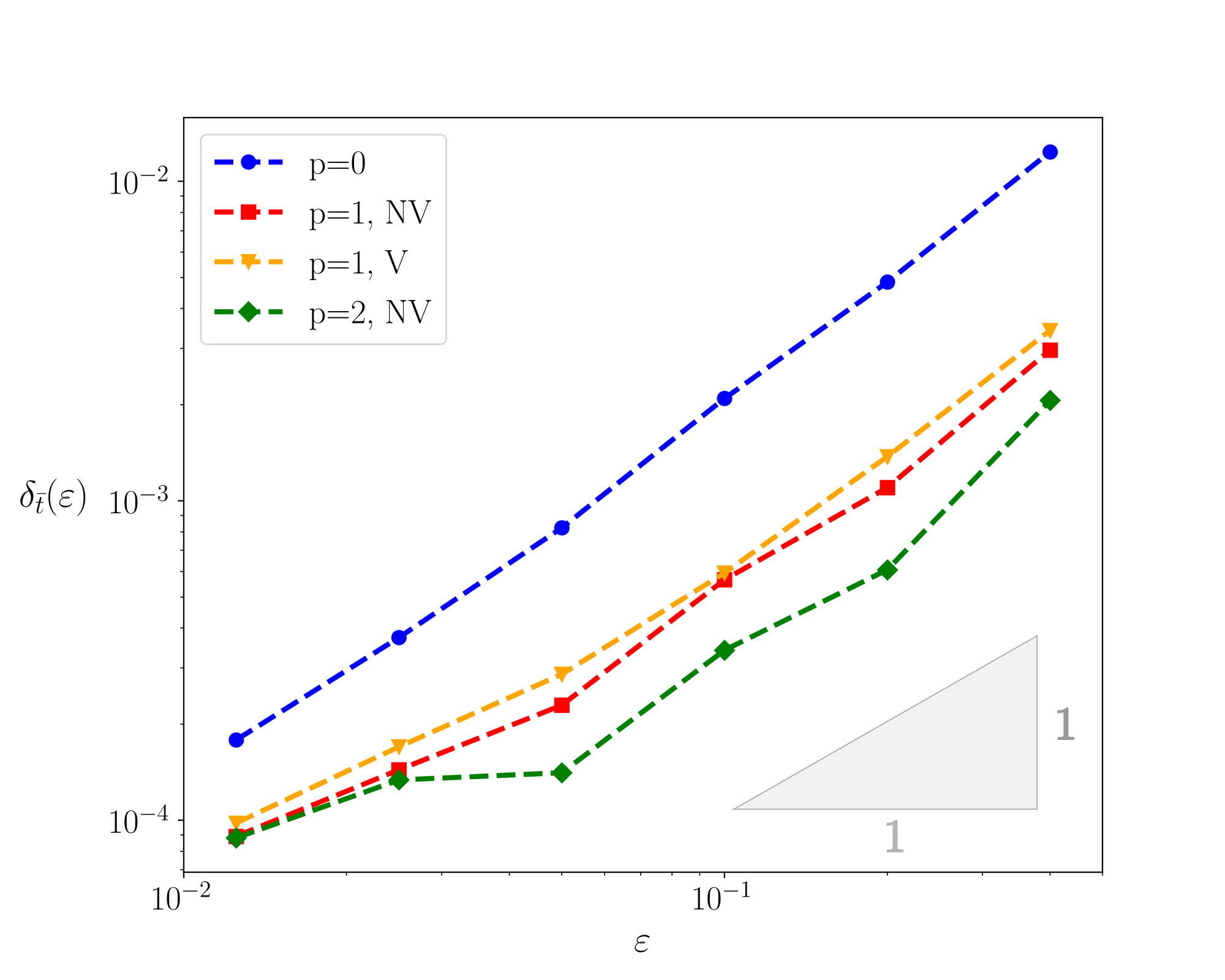} 
\caption{Convergence for $\delta_{\bar{t}} (\varepsilon)$ at $\bar{t}=0.0003$ using DCH model ($p = 0$) and DDCH models ($p = 1$ Model V, Model NV and $p = 2$ Model NV).}
\label{fig:figureConv}
\end{figure}

The convergence rate is considered in Figure~\ref{fig:figureConv} using the relative deviation from the SI solution at time $\bar{t}$, as
\begin{equation}
    \delta_{\bar{t}}(\varepsilon)= \frac{|A_x(\bar{t},\varepsilon)^{\rm (D)DCH}-A_x(\bar{t})^{\rm SI}|}{A_x(\bar{t})^{\rm SI}}.
\end{equation}
This quantity is shown for an early stage of the simulation focusing on the fast dynamics. For all DDCH models $p = 0$, $p = 1$ and $p = 2$ we obtain linear convergence. However, with a significantly smaller error for $p = 1$ and $p = 2$. This confirms previous results \cite{Backofen2019} and extends them to Model V. In any case, using the DDCH models allows to consider $\varepsilon$-values to be at least doubled if compared with the classical DCH model, to reach the same accuracy.

\begin{figure}[h]
    \center
\includegraphics*[width = 1.0 \textwidth ]{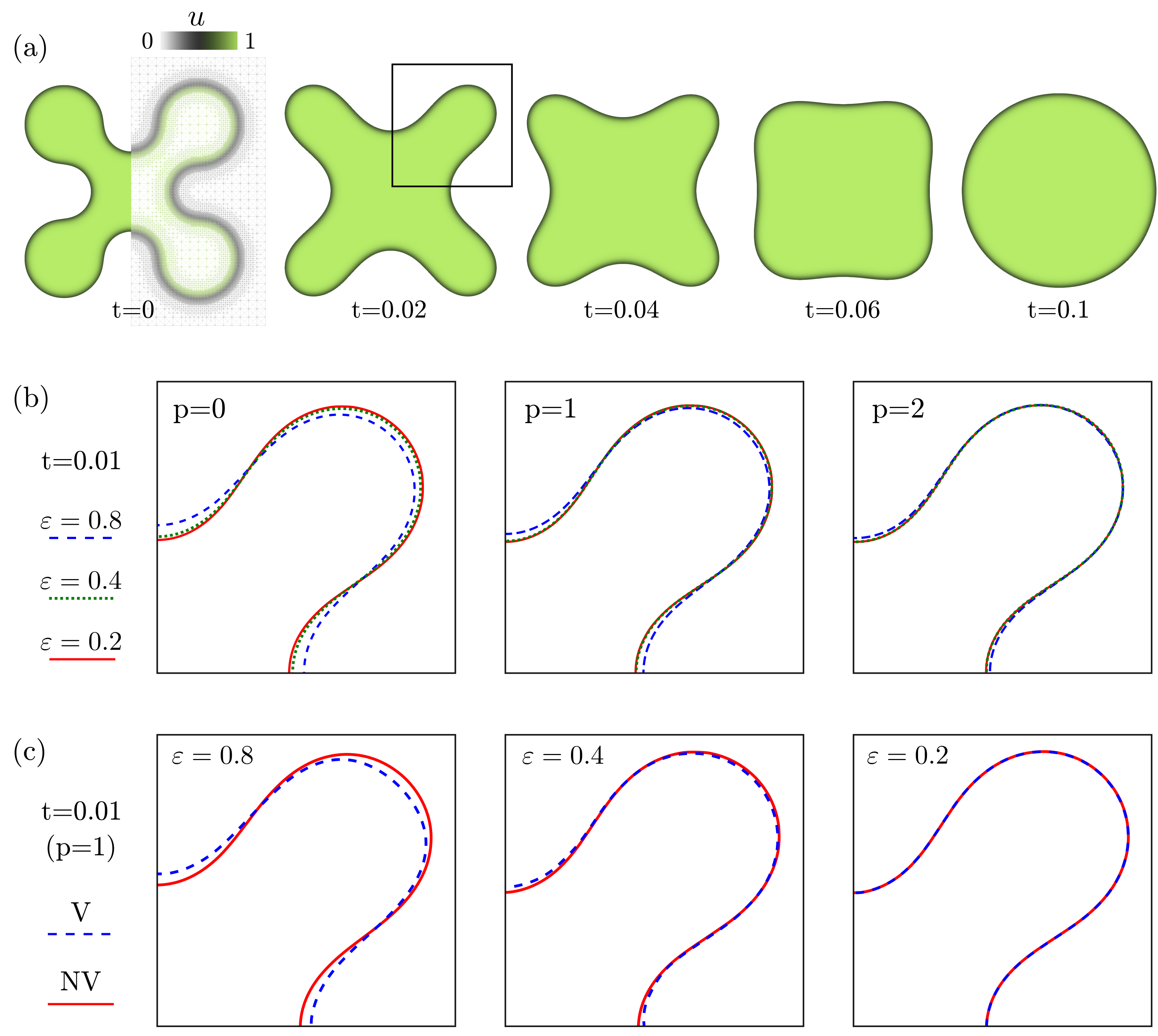} 
\caption{Evolution of a four-fold shape with positive and negative curvature regions. (a) Initial stage, reporting the region where $u>0.5$ (left) and the mesh adopted for the simulation (right) at $t=0$, with four representative stages to  illustrate the dynamics. A stationary shape is obtained (approximately) for $t>0.1$. (b) Comparison of the surface profile $u\sim 0.5$ at $t=0.01$ obtained by Model NV for different values of $\varepsilon$ and $p$. (c) Comparison of the surface profile $u\sim 0.5$ at time $t=0.01$ obtained by Models V and NV with $p=1$, for different values of $\varepsilon$. Panels (b) and (c) show the northeast region corresponding to the black square of panel (a) (see $t=0.02$).}
    \label{fig:figure5}
    \end{figure}

The behavior of the DDCH models discussed so far holds true also when considering more complex surface profiles. The case of the evolution by surface diffusion of a profile where also the sign of the local curvature changes is illustrated in Figure~\ref{fig:figure5}.\footnote{The profile of Figure~\ref{fig:figure5} is known as \textit{Camunian Rose}, a stylized version of iron-age rock carvings found in northern Italy~\cite{Farina1997}.} Figure~\ref{fig:figure5}(a) at $t=0$ shows the considered profile obtained with $\varepsilon=0.2$, along with the computational mesh adopted for the numerical simulations. In the same panel, four representative stages during the evolution are also shown. The shapes in Figure~\ref{fig:figure5}(a) correspond to the region $u>0.5$ with a color map highlighting the extension of the interface region. Figure~\ref{fig:figure5}(b) shows the different surface profiles, namely the $u = 0.5$ level curve, obtained with different values of $\varepsilon$ and different choice of $p$ for the Model NV. The increase of accuracy with increasing $p$ is qualitatively confirmed here. Figure~\ref{fig:figure5}(c) shows that Model V ($p=1$) delivers very similar results to Model NV as expected from the convergence results illustrated in Figure~\ref{fig:figureConv}.

    \section{Discussion and Conclusions}
    \label{sec:conclusion}
    
In this paper we have introduced a new variational model, Model V, describing surface diffusion. This model has an energy, though this energy is defined via a singular restriction function that must be chosen carefully. One of the main purposes for the creation of the new model, was to connect it, via certain  approximations, to a well known, and well used, non-variational model, Model NV. Both Models V and NV are shown to converge to sharp interface surface diffusion, via the formal method of matched asymptotics. This convergence is further supported by our numerical results, which suggest that the convergence rates with respect to $\varepsilon$ are all order 1, that is, the convergence rates are all linear in $\varepsilon$.

Both models are doubly degenerate Cahn-Hilliard (DDCH) equations, and there are relative strengths and weaknesses associated to each. The new model, Model V, is more complicated to solve numerically, but has an associated energy, which makes its connection to other physics though additional energy terms, more seamless.  Like for Model NV, the solutions to Model V, at a fixed $\varepsilon$ value and $p$ value,  more accurately approximate  solutions of the sharp interface model, compared to the singly degenerate Cahn-Hillard (DCH) model (where, the only degeneracy is associated to the mobility). Solutions to Model NV, it seems, are slightly more accurate approximations than those of Model V for the $p=1$ case; and solutions obtained from Model NV with $p=2$ are more accurate than any of the solutions computed from Model V. 

Several open questions related to both models remain, questions related to existence, uniqueness, and regularity of solutions, as well as to the positivity preserving property. In our computations -- where we used the regularized versions of the equations, since our schemes could not maintain positivity --  solutions to Model V were closer to being positive than solutions to Model NV for the $p=1$ case (data not shown). In other words, the overshoots of the values outside of the physically realistic range, $0 \le u \le 1$, were smaller in a point-wise sense for Model V. Indeed, Model V has a singular energetic mechanism for positivity that Model NV may not possess. However, either way of including a restriction function seems to help in maintaining positivity, when compared to the standard singly degenerate Cahn-Hilliard model for surface diffusion.  Our numerical experiments suggest that solutions to both models will remain positive, in the sense that reducing the regularization reduces the point-wise overshoot of the values outside of the physically realistic range, $0 \le u \le 1$. In future works, we plan to report on energy stable and positivity preserving numerical schemes for Model V, as well as some other theoretical issues, including the existence of weak solutions. Moreover, for the extension of Model V to the anisotropic setting, crucial to account for faceting of surfaces evolving by surface diffusion we refer to \citen{SalvalaglioDDCHaniso}.
    
    \section*{Acknowledgments}
 
This research was partially funded by the EU H2020 FET-OPEN project microSPIRE (ID: 766955) and by the EU H2020 FET-OPEN project NARCISO (ID: 828890). We gratefully acknowledge the computing time granted by the John von Neumann Institute for Computing (NIC) and provided on the supercomputer JURECA at J\"ulich Supercomputing Centre (JSC), within the Project no. HDR06, and by the Information Services and High Performance Computing (ZIH) at the Technische Universit\"at Dresden (TUD). SMW acknowledges generous financial support from the US National Science Foundation, through grants NSF-DMS 1719854 and NSF-DMS 2012634, and acknowledges some enlightening discussions with Prof.~Shibin Dai  on the topics of degenerate mobilities in the Cahn-Hilliard equation and with Prof.~Cheng Wang on the topic of energy stable and positivity preserving methods for Cahn-Hilliard-deGennes-type models. Some of this work was performed while the third author (SMW) was supported as a Dresden Senior Fellow at TUD for three months in 2017 in the group of the second author (AV). SMW gratefully acknowledges this funding and the hospitality of AV and TUD during his visit.

    \section*{Conflicts of Interest Statement}
The authors certify that they have no affiliations with or involvement in any organization or entity with any financial or non-financial interest in the subject matter discussed in this manuscript.

    \appendix

	\section{Matched Asymptotic Analysis of Model V}
	\label{sec-asymptotics-model-v}
	
In this appendix, we will prove, using the method of matched asymptotic expansions, that solutions of the variational DDCH system \eqref{eqn-ch-little-g-1} -- \eqref{eqn-ch-little-g-2}, which we refer to as Model V, converge to the solutions of the sharp interface model of surface diffusion \eqref{eqn-surf-diff}. To keep the discussion simple, we consider only the case  for which space is two dimensional. In this setting, the surface diffusion of a one-dimensional evolving curve is modelled by the equation
    \[
v = \partial_{ss}\kappa,
    \]
where $s$ is the arc length of $\Sigma$, and $\kappa$ is the curvature.  Herein we assume that $0\le p < 2$ in the definition of the restriction function $g_\alpha$ and we set the regularization to zero ($\alpha = 0$). We closely follow the analyses in \citen{cahn1996,Ratz2006}.

	\subsection{Interface Coordinates and Asymptotic Expansions}
	
We suppose that $\Omega\subset\mathbb{R}^2$ is a bounded open set, for example, a rectangular domain. Let us assume that there is a simply-connected domain $\Omega^0_0\subset \Omega$ that has a smooth boundary $\Sigma_0 := \partial\Omega_0$. Assume that $u_{0,\varepsilon}:\Omega \to [0,1]$ is a smooth function with the property that 
    \[
\Sigma_0 = \left\{ \bfx \ \middle| \ u_{0,\varepsilon}(\bfx) = 0.5\right\} \quad \mbox{and} \quad \Omega^0_0 = \left\{ \bfx \ \middle| \ u_{0,\varepsilon}(\bfx) < 0.5\right\} .
    \]
The idea is that, as $\varepsilon\searrow 0$, $u_{0,\varepsilon}$ converges to the characteristic function $1_{\Omega\setminus \Omega_0}$. We assume that $u$ is the solution to Model V, that is, \eqref{eqn-ch-little-g-1} -- \eqref{eqn-ch-little-g-2} with initial data $u_{0,\varepsilon}$ and 
    \[
\Sigma(t,\varepsilon) = \left\{ \bfx \ \middle| \ u(\bfx,t) = 0.5\right\} \quad \mbox{and} \quad \Omega^0(t,\varepsilon) = \left\{ \bfx \ \middle| \ u(\bfx,t) < 0.5\right\} .
    \]
Clearly $\Sigma_0 = \Sigma(0;\varepsilon)$. We assume that $\Sigma(t,\varepsilon)$ remains smooth and there is always a neighborhood $N(t;\varepsilon) \supset \Sigma(t;\varepsilon)$ such that there is a well-defined signed distance function $r(x,y,t;\varepsilon)$. There is a well-definite interior, where $r<0$ and corresponding to $\Omega^0(t,\varepsilon)$, and an exterior, where $r>0$ and corresponding to $\Omega \setminus \Omega^0(t,\varepsilon)$.

We assume for simplicity that $\Sigma(t,\varepsilon)$ is closed, and its length, denoted $L(t;\varepsilon)$, is finite. Thus $\Sigma(t,\varepsilon)$ can be parameterized by its arc length, $s\in [0, L(t;\varepsilon)]$. The parameterization is denoted by the vector function
	\[
\bfc = \bfc( \, \cdot \, , t; \varepsilon) :  \mathbb{R} \longrightarrow 	\Sigma(t,\varepsilon),
	\]
which is an $L(t;\varepsilon)$-periodic, twice continuously differentiable function that is a bijection when its domain is restricted to any half-open interval $[a,b)$ of length $L(t;\varepsilon)$.

On the interface $\Sigma(t;\varepsilon)$, there are well-defined unit tangent and unit normal vector functions, which we denote by $\bftau(s,t;\varepsilon)$ and $\bfn(s,t;\varepsilon)$. The unit tangent vector points in the direction of increasing arc length, while the unit normal points in the direction of increasing interface distance $r$. In fact, in a small enough neighborhood 
	\[
N_\rho(t;\varepsilon) := \left\{(x,y)\in\Omega \ \middle| \ |r(x,y,t;\varepsilon)| < \rho  \right\} \subseteq  N(t;\varepsilon),
	\]
$s$ and $r$ are local coordinates near the interface $\Sigma(t;\varepsilon)$. The set $N_\rho(t;\varepsilon)$ is called the \textit{narrow band neighborhood of the interface}, or just the \textit{narrow band}, for short. In particular, for every $\bfx = (x,y) \in N_\rho(t;\varepsilon)$, there are points $r\in (-\rho,\rho)$ and $s\in [0,L(t;\varepsilon))$ such that
	\[
(x,y) = \bfx(r,s) = \bfc(s,t;\varepsilon) + r(x,y,t;\varepsilon)\bfn(s,t;\varepsilon).	
	\] 
In other words, for every $(x,y) \in N_\rho(t;\varepsilon)$, there is a unique pair 
	\[
(r,s) \in \	B(t; \varepsilon) :=  \left\{ (r,s) \ \middle| \  -\rho < r < \rho , \ 0\le s < L(t;\varepsilon) \right\}.
	\]
We can now define transformed dependent variables in the narrow band via composition: for all $(r,s) \in	B(t; \varepsilon)$, define
	\begin{align*}
\hat{u}(r,s,t;\varepsilon) & := u(\bfx(r,s),t;\varepsilon),
	\\
\hat{w}(r,s,t;\varepsilon) & := w(\bfx(r,s),t;\varepsilon) .
	\end{align*}
	
To facilitate our asymptotic analysis, we introduce the stretched interface distance variable
	\[
z = \frac{r}{\varepsilon} .
	\]
Subject to this transformation, we define, for all $(r,s) \in B(t; \varepsilon)$,  the stretched variables
	\begin{align*}
U(z,s,t;\varepsilon) & := \hat{u}(r,s,t;\varepsilon) ,
	\\
W(z,s,t;\varepsilon) & := \hat{w}(r,s,t;\varepsilon) .
	\end{align*}
We assume that the following Taylor expansions are valid
	\begin{align}
u(x,y,t;\varepsilon) &= u_0(x,y,t) + \varepsilon u_1(x,y,t) + \varepsilon^2 u_2(x,y,t) + \cdots ,
	\label{u-outer-x-y}
	\\
\hat{u}(r,s,t;\varepsilon) &= \hat{u}_0(r,s,t) + \varepsilon \hat{u}_1(r,s,t) + \varepsilon^2 \hat{u}_2(r,s,t) + \cdots,
	\label{u-outer-r-s}
	\\
U(z,s,t;\varepsilon) &= U_0(z,s,t) + \varepsilon U_1(z,s,t) + \varepsilon^2 U_2(z,s,t) + \cdots,
	\label{u-inner-z-s}
	\\
w(x,y,t;\varepsilon) &= w_0(x,y,t) + \varepsilon w_1(x,y,t) + \varepsilon^2 w_2(x,y,t) + \cdots ,
	\label{w-outer-x-y}
	\\
\hat{w}(r,s,t;\varepsilon) &= \hat{w}_0(r,s,t) + \varepsilon \hat{w}_1(r,s,t) + \varepsilon^2 \hat{w}_2(r,s,t) + \cdots,
	\label{w-outer-r-s}
	\\
W(z,s,t;\varepsilon) &= W_0(z,s,t) + \varepsilon W_1(z,s,t) + \varepsilon^2 W_2(z,s,t) + \cdots.
	\label{w-inner-z-s}
	\end{align}
Expansions \eqref{u-outer-x-y}, \eqref{u-outer-r-s}, \eqref{w-outer-x-y}, and \eqref{w-outer-r-s} are called outer expansions; while expansions \eqref{u-inner-z-s} and \eqref{w-inner-z-s} are called inner expansions. The following \textit{matching conditions} are employed to join the two types of expansions together:
	\begin{align}
\lim_{r\to \pm 0}\hat{u}(r,s,t;\varepsilon) = \lim_{z\to\pm \infty} U(z,s,t;\varepsilon),
	\\
\lim_{r\to \pm 0}\hat{w}(r,s,t;\varepsilon) = \lim_{z\to\pm \infty} W(z,s,t;\varepsilon).
	\end{align}

We denote by $v$ the normal velocity of the $\Sigma(t,\varepsilon)$ and, by $\kappa$, the curvature of $\Sigma(t,\varepsilon)$. Then, we have the following well-known relations in two dimensions:
	\[
\partial_s \bftau = -\kappa \bfn \quad \mbox{and} \quad \partial_s \bfn = \kappa\bftau.	
	\]
As a consequence, using the chain rule,
	\begin{align*}
\partial_s \hat{u} & = \nabla_{\bfx} u \cdot \frac{\partial\bf x}{\partial s} 
	 = \nabla_{\bfx} u \cdot\left(\partial_s \bfc +r\partial_s\bfn \right) 
	 = \nabla_{\bfx} u \cdot\left(\partial_s \bfc +r\kappa\bftau  \right) .
	\end{align*}
Since we are parameterizing by arc-length, it follows that 
	\[
\partial_s \bfc = \bftau  ,	
	\]
and, therefore, 
	\[
\partial_s \hat{u} = \left(1 +r\kappa  \right) \nabla_{\bfx} u \cdot \bftau.
	\]
Equivalently,
	\[
\nabla_{\bfx} u \cdot \bftau = \frac{\partial_s \hat{u}}{1 +r\kappa }. 
	\]
In a similar way, we can derive the identity
	\[
\nabla_{\bfx} u \cdot \bfn = \partial_r \hat{u}.
	\]

Using these and related identities, we can express the standard operators in the transformed coordinate system $(z,s)$:
	\begin{align*}
\partial_t u & = \partial_t\hat{u} -v \partial_r\hat{u} -(1+r\kappa)^{-1}\partial_s\hat{u}\left(\partial_t\bfc +r\partial_t\bfn  \right)\cdot \bftau
	\\
& = \partial_t U -\varepsilon^{-1} v \partial_z U -\left(1+\varepsilon z \kappa\right)^{-1}\partial_s U \left(\partial_t\bfc +\varepsilon z \partial_t\bfn  \right) \cdot \bftau,
	\\
\nabla u & = \partial_r \hat{u} \bfn + (1+r\kappa)^{-1} \partial_s \hat{u} \bftau
	\\
& = \varepsilon^{-1} \partial_z U \bfn + (1+\varepsilon z \kappa)^{-1} \partial_s U \bftau,
	\\
\nabla\cdot \bff & = \partial_r(\hat{\bff}\cdot \bfn) +(1+r\kappa)^{-1}\left(\partial_s(\hat{\bff}\cdot\bftau) +\kappa\hat{\bff}\cdot \bfn  \right)
	\\
& = \varepsilon^{-1}\partial_z(\bfF\cdot \bfn) +(1+\varepsilon z \kappa)^{-1}\left(\partial_s(\bfF\cdot\bftau) +\kappa\bfF\cdot \bfn  \right) ,
	\end{align*} 
where $\bff$ is any vector function that transforms as $\hat\bff(r,s,t) = \bff(\bfx(r,s),t)$ and $\bfF(z,s,t) = \hat\bff(r,s,t)$. The full derivations of the equations above can be found in~\citen{Ratz2006}. Combining the expansions above, we can express the Laplacian operator as
	\begin{align*}
\Delta u & = \varepsilon^{-2} \partial_{zz}U +\varepsilon^{-1}(1+\varepsilon z \kappa)^{-1} \kappa \partial_z U -\varepsilon z \partial_s\kappa(1+\varepsilon z\kappa )^{-3} \partial_s U+(1+\varepsilon z \kappa)^{-2} \partial_{ss}U .
	\end{align*} 
Likewise,
	\begin{align*}
|\nabla u|^2 = \nabla u\cdot \nabla u & = \left(  \varepsilon^{-1} \partial_z U \bfn + (1+\varepsilon z \kappa)^{-1} \partial_s U \bftau \right)\cdot \left(\varepsilon^{-1} \partial_z U \bfn + (1+\varepsilon z \kappa)^{-1} \partial_s U \bftau\right)
	\\
& =  \varepsilon^{-2} \left( \partial_z U \right)^2 + (1+\varepsilon z \kappa)^{-2} \left(\partial_s U\right)^2.
	\end{align*}
	
	\subsection{Outer Expansion}
The evolution equations become singular in the outer region. While the energy for the pure states can be defined in reasonable way, in particular,
    \[
F[u\equiv 0] = 0 = F[u\equiv 1] ,
    \]
the gradient equations are undefined for the pure states. This is analogous to the deep-quench limit examined in \citen{cahn1996}. An expansion of solutions may not be available in the usual sense. Without more information, we will assume the following reasonable form of the solution in the outer regions:
	\[
u = 0 + \varepsilon \beta + \mathcal{O}(\varepsilon^2) \ge 0  \quad \mbox{and} \quad u = 1 - \varepsilon \beta + \mathcal{O}(\varepsilon^2) \le 1,
	\]
where $\beta\ge 0$ is a constant. This is enough for our analysis to go through, though we leave open the possibility that $\beta = 0$, as was assumed in \citen{cahn1996} for the deep-quench limit case. The structure of the outer solutions corresponds to matching conditions for the inner expansions of the form
	\[
\lim_{z\to-\infty} U_0 = 0, \quad 	\lim_{z\to\infty} U_0 = 1, \quad  \lim_{z\to-\infty} U_1 = \beta, \quad \mbox{and} \quad \lim_{z\to\infty} U_1 = -\beta.
	\]
	
	\subsection{Inner Expansion}
	
	\subsubsection{Equation~\eqref{eqn-ch-little-g-2} at $\mathcal{O}(\varepsilon^{-1})$}
	
We will assume that in the inner region, the solution stays away from the pure states:
	\[
0  <  u < 1.	
	\]
Let us make an expansion for Eq.~\eqref{eqn-ch-little-g-2}, which is rewritten slightly as
	\[
w  =  g'(u) \left(\frac{1}{\varepsilon} f(u) - \frac{\varepsilon}{2}|\nabla u |^2\right)  + g(u) \left(\frac{1}{\varepsilon} f'(u)  - \varepsilon  \Delta u\right) ,
	\]
where for simplicity, we have dropped the subscript $\alpha$ on $g_\alpha$.  Subscripts will have a different meaning in this section, as we will see. We have
	\begin{align}
W_0 + \mathcal{O}(\varepsilon) & = \left(g'_0 + \varepsilon g'_1 + \mathcal{O}(\varepsilon^2) \right) \left(\varepsilon^{-1}f_0 + \varepsilon^0 f_1 - \varepsilon^{-1}\frac{1}{2}(\partial_zU_0)^2 -\varepsilon^0 \partial_zU_0\partial_z U_1 + \mathcal{O}(\varepsilon) \right)
	\nonumber
	\\
& \quad + \left(g_0 + \varepsilon g_1 + \mathcal{O}(\varepsilon^2) \right)\left( \varepsilon^{-1}f'_0 + \varepsilon^0 f'_1 - \varepsilon^{-1}\partial_{zz}U_0 - \varepsilon^0\partial_{zz}U_1 - \varepsilon^0\kappa\partial_z U_0 + \mathcal{O}(\varepsilon) \right)
	\nonumber
	\\
& = \left(g'_0 + \varepsilon g'_1 + \mathcal{O}(\varepsilon^2) \right) \left(\varepsilon^{-1} \left( f_0 - \frac{1}{2}(\partial_zU_0)^2\right)+ \varepsilon^0\left( f_1  - \partial_zU_0\partial_z U_1\right) + \mathcal{O}(\varepsilon) \right)
	\nonumber
	\\
& \quad + \left(g_0 + \varepsilon g_1 + \mathcal{O}(\varepsilon^2) \right)\left( \varepsilon^{-1}\left(f'_0- \partial_{zz}U_0\right) + \varepsilon^0\left( f'_1  - \partial_{zz}U_1 -\kappa\partial_z U_0\right) + \mathcal{O}(\varepsilon) \right).
	\label{inner-ch-2}
	\end{align}
Here
	\[
g_0 := g(U_0), \quad g_1 = g'(U_0) U_1 , \quad 	f_0 := f(U_0), \quad f_1 = f'(U_0) U_1
	\]
and
	\[
g'_0 := g'(U_0), \quad g'_1 = g''(U_0) U_1 , \quad 	f'_0 := f'(U_0), \quad f'_1 = f''(U_0) U_1
	\]
	
At $\mathcal{O}(\varepsilon^{-1})$ we get
	\[
0 = g'(U_0) \left( f_0 - \frac{1}{2}(\partial_zU_0)^2\right) + g(U_0)\left(f'_0- \partial_{zz}U_0\right).
	\]
Multiplying by $\partial_z U_0$, we get a perfect derivative:
	\begin{align*}
0 & = g'(U_0) \partial_z U_0 \left( f_0 - \frac{1}{2}(\partial_zU_0)^2\right) + g(U_0)\left(f'_0- \partial_{zz}U_0\right)\partial_z U_0 
	\\
& = \partial_z\left[ g(U_0) \left(f_0 - \frac{1}{2}(\partial_zU_0)^2 \right) \right].
	\end{align*}
Therefore
	\[
g(U_0) \left(f(U_0) - \frac{1}{2}(\partial_zU_0)^2 \right) = C_1(s,t) .
	\]
Using the matching conditions, 
	\[
\lim_{z\to -\infty} U_0 = 0 \quad \mbox{and} \quad 	\lim_{z\to \infty} U_0 = 1,
	\]
and the fact that $0\le p < 2$, in which case $f(U_0)$ goes to zero faster than $g(U_0)$ goes to infinity, we conclude that $C_1$ must be zero. Since $g(U_0) \ne 0$, it must be that 
	\[
f(U_0) - \frac{1}{2}(\partial_zU_0)^2 = 0 = f'(U_0) -\partial_{zz} U_0 .
	\]
The solution we seek is
	\[
U_0 = \frac{1}{2}\left(1+\tanh\left(3z \right) \right)	= \frac{1}{2}\left(1+\tanh\left(\frac{3r}{\varepsilon} \right) \right) .
	\]
Now, using substitution, with
    \[
dU_0 = 6 U_0 (1-U_0) dz ,    
    \]
and the matching conditions, we can compute a couple of common integrals that we will meet later: first
	\begin{align}
\int_{-\infty}^\infty \left(\partial_z U_0\right)^2 d z & = \int_{-\infty}^\infty 36U_0^2(1-U_0)^2 dz = 6 \int_0^1 U_0(1-U_0) \, dU_0 = 1.
    \label{integral-U0-1}
	\end{align}
Similarly,
    \begin{equation}
\int_{-\infty}^\infty \partial_z U_0 \, d z = \int_{-\infty}^\infty 6 U_0(1-U_0) dz = \int_0^1 \, dU_0 = 1.
    \label{integral-U0-2}
    \end{equation}

	\subsubsection{Equation~\eqref{eqn-ch-little-g-1} at $\mathcal{O}(\varepsilon^{-3})$}
	
Now, let us consider Eq.~\eqref{eqn-ch-little-g-1}. Inserting the inner expansion, we get
	\begin{align}
\partial_t U & - \varepsilon^{-1} v \partial_z U -\left(1+\varepsilon z \kappa\right)^{-1}\partial_s U \left(\partial_t\bfc +\varepsilon z \partial_t\bfn  \right) \cdot \bftau
	\nonumber
	\\
& = \varepsilon^{-3}\partial_z(M(U)\partial_z W  ) +\varepsilon^{-1}(1+\varepsilon z \kappa)^{-1}\left[\partial_s\left\{M(U)(1+\varepsilon z\kappa)\partial_s W \right\} +\kappa M(U)\varepsilon^{-1}\partial_z W   \right].
	\label{eq-ch-1-inner}
	\end{align}	
Observe the following expansion
	\[
M(U) = M_0 +\varepsilon M_1 +\varepsilon^2 M_2 + \varepsilon^3 M_3 + \mathcal{O}(\varepsilon^4)	,
	\]
where
	\[
M_0  = M(U_0), \quad M_1  = M'(U_0) U_1 , \quad M_2  = \frac{1}{2}\left(M''(U_0)U_1^2 + 2 M'(U_0) U_2\right), \quad \cdots	
	\]
Therefore, in \eqref{eq-ch-1-inner} at $\mathcal{O}(\varepsilon^{-3})$, we have
	\[
0 = \partial_z(M_0\partial_z W_0  ).
	\]
This implies that
	\[
M(U_0)\partial_z W_0 = C_2(s,t).
	\]
Using the matching condition, it is clear that
	\[
\lim_{z\to \pm \infty} M(U_0(z,s,t)) = 0,	
	\]
which implies that $C_2 = 0$. Therefore,
	\[
\partial_z W_0 = 0.	
	\]
This implies that
	\[
W_0 = C_3(s,t) .	
	\]
	
	\subsubsection{Equation~\eqref{eqn-ch-little-g-1} at $\mathcal{O}(\varepsilon^{-2})$}
	
Using \eqref{eq-ch-1-inner}, at $\mathcal{O}(\varepsilon^{-2})$, we have
	\[
0 = \partial_z(M_0\partial_z W_1 + M_1\partial_z W_0 + \kappa M_0\partial_z W_0 ) = \partial_z(M_0\partial_z W_1) .
	\]
As above, we have
	\[
\partial_z W_1 = 0,	
	\]
which implies that
	\[
W_1 = C_4(s,t) .	
	\]

	\subsubsection{Equation~\eqref{eqn-ch-little-g-1} at $\mathcal{O}(\varepsilon^{-1})$}
	
Finally, in \eqref{eq-ch-1-inner} at $\mathcal{O}(\varepsilon^{-1})$, we have, after using $\partial_z W_0 = 0 =\partial_z W_1$,
	\begin{align*}
-v \partial_z U_0 = \partial_z(M_0 \partial_z W_2) + \partial_s(M_0 \partial_s W_0) .
	\end{align*}
Integrating along the $z$-axis, 
	\begin{equation}
-v  = \int_{-\infty}^\infty \partial_z(M_0 \partial_z W_2) dz + \partial_s\left(\int_{-\infty}^\infty M(U_0)dz \partial_s W_0\right)  = \partial_{ss} W_0 ,
	\label{velocity-equation}
	\end{equation}
where we have used \eqref{integral-U0-2} and the fact that
	\[
\int_{-\infty}^\infty M(U_0)dz = \int_{-\infty}^\infty 36U_0^2(1-U_0)^2 dz =  1,	
	\]
which can be deduced from \eqref{integral-U0-1}.
	
	\subsubsection{Equation~\eqref{eqn-ch-little-g-2} at $\mathcal{O}(\varepsilon^{0})$}

Now, going back to Eq.~\ref{eqn-ch-little-g-2}, using \eqref{inner-ch-2}, at order $\mathcal{O}(\varepsilon^{0})$, we have
	\begin{align*}
W_0 & = g'_0 \left( f_1  - \partial_zU_0\partial_z U_1\right) +g'_1\left( f_0 - \frac{1}{2}(\partial_zU_0)^2\right) 
	\\
& \quad + g_0 \left( f'_1  - \partial_{zz}U_1 -\kappa\partial_z U_0\right) + g_1 \left(f'_0- \partial_{zz}U_0\right).
	\end{align*}
This simplifies to 
	\[
W_0 + g(U_0) \kappa\partial_z U_0  = 	g'(U_0) \left( f'(U_0)U_1  - \partial_zU_0\partial_z U_1\right) + g(U_0) \left( f''(U_0) U_1  - \partial_{zz}U_1 \right).
	\]
Multiplying by $\partial_zU_0$, integrating, and being careful about boundary terms, we have
	\begin{align*}
\int_{-\infty}^\infty	\left( W_0 + g(U_0) \kappa\partial_z U_0\right)\partial_zU_0 \, dz & = \int_{-\infty}^\infty g'(U_0) f'(U_0)U_1 \partial_z U_0\, dz 
	\\
& \quad - \int_{-\infty}^\infty g'(U_0)  \partial_z U_0 \partial_z U_1 \partial_z U_0\, dz
	\\
& \quad + \int_{-\infty}^\infty g(U_0) f''(U_0) U_1 \partial_z U_0 \, dz 
	\\
& \quad - \int_{-\infty}^\infty g(U_0) \partial_{zz} U_1 \partial_z U_0 \, dz 
	\\
& =  \int_{-\infty}^\infty \partial_z \left[g(U_0) f'(U_0)\right] U_1 \, dz
	\\
& \quad - \int_{-\infty}^\infty \partial_z \left[g(U_0) \partial_z U_1  \right] \partial_z U_0 \, d z
	\\
& = -\int_{-\infty}^\infty g(U_0) f'(U_0)  \partial_zU_1 \, dz 
	\\
& \quad + \int_{-\infty}^\infty g(U_0) \partial_z U_1  \partial_{zz} U_0 \, d z  +B
	\\
& = - \int_{-\infty}^\infty g(U_0) \left(f'(U_0) - \partial_{zz} U_0  \right) \partial_z U_1 \, dz + B 
    \\
& =  B,
	\end{align*}
where
	\[
B  :=  \left[g(U_0)\left(f'(U_0) U_1 - \partial_z U_1\partial_z U_0\right)\right]_{z=-\infty}^{z=\infty} = -2\beta g(0) f'(0) = 12 \beta,
	\]
using the matching conditions for $U_0$ and $U_1$. The important point is that $B$ is a constant; in fact, it may be zero, but we only need it to be constant. Observe that
    \begin{equation}
\int_{-\infty}^\infty g(U_0) \left(\partial_z U_0\right)^2 \, dz = \frac{6}{\gamma} \int_0^1U_0^{1-p} (1-U_0)^{1-p}\, dU_0 = 1  ,
    \label{integral-U0-3}
    \end{equation}
since $0\le p<2$ and we use the energy  normalization \eqref{eqn-energy-normalize-gamma}. Using \eqref{integral-U0-2} and \eqref{integral-U0-3}, and the fact that $W_0$ and $\kappa$ are independent of $z$, it follows that 
    \[
\int_{-\infty}^\infty	\left( W_0 + g(U_0) \kappa\partial_z U_0\right)\partial_zU_0 \, dz = W_0 + \kappa,
    \]
and, therefore, 
	\begin{equation}
W_0 = - \kappa + B ,
	\label{chemical-pot-equation}
	\end{equation}
where $B$ is a constant. Combining Equations \eqref{velocity-equation} and \eqref{chemical-pot-equation}, we have the desired result
	\[
v = \partial_{ss}\kappa .	
	\]
	
	\section{Asymptotic Analysis of Model NV}
	\label{sec-asymptotics-model-nv}

The analysis for Model NV is similar to that of Model V, and we use analogous notation. All of the details for the case $p=2$ can be found in \citen{Ratz2006}, and we only briefly summarize our results.  We can prove that solutions converge, as $\varepsilon\searrow 0$, to solutions of the SI problem $v = \partial_{ss} \kappa$
provided that in the inner expansion 
    \[
 \int_{-\infty}^\infty G( U_0)\partial_z  U_0 dz = 1 \quad \mbox{and} \quad \int_{-\infty}^\infty \left(\partial_z  U_0 \right)^2 dz = 1.  
    \]
To see this, consider the inner expansion for Eq.~\eqref{eqn-ch-G-2}, which is rewritten slightly as
	\[
G( u) w  = \frac{1}{\varepsilon} f'( u)  - \varepsilon  \Delta  u  ,
	\]
where for simplicity, we have dropped the subscript 0 on $G_0$, as above.  Subscripts will again have a different meaning. We have
	\begin{align}
G( U_0)W_0 + \mathcal{O}(\varepsilon) & =  \varepsilon^{-1}\left(f'_0- \partial_{zz} U_0\right) + \varepsilon^0\left( f'_1  - \partial_{zz} U_1 -\kappa\partial_z  U_0\right) + \mathcal{O}(\varepsilon) ,
	\label{inner-ch-G-2}
	\end{align}
where
	\[
G_0 := G( U_0), \quad G_1 = G'( U_0)  U_1 , \quad f'_0 := f'( U_0), \quad f'_1 = f''( U_0)  U_1 .
	\]
At $\mathcal{O}(\varepsilon^{-1})$ we get
	\[
0 =  f'_0- \partial_{zz} U_0,
	\]
which guarantees that 
	\[
 U_0 = \frac{1}{2}\left(1+\tanh\left(3z \right) \right) \quad \mbox{and} \quad \int_{-\infty}^\infty \left(\partial_z  U_0 \right)^2 dz = 1.
	\]
At order $\mathcal{O}(\varepsilon^{0})$, we have
	\begin{align*}
G( U_0)W_0 & =  f'_1  - \partial_{zz} U_1 -\kappa\partial_z  U_0 .
	\end{align*}
This simplifies to 
	\[
G( U_0)W_0 +  \kappa\partial_z  U_0  =  f''( U_0)  U_1  - \partial_{zz} U_1 .
	\]
Multiplying by $\partial_z U_0$ and integrating, we get
    \[
\int_{-\infty}^\infty G( U_0)W_0 \partial_z U_0 \ dz + \int_{-\infty}^\infty \kappa\left(\partial_z  U_0  \right)^2 dz =  \int_{-\infty}^\infty \left(f''( U_0)  U_1  - \partial_{zz} U_1\right)\partial_z U_0 \, dz = 0.
    \]

In \citen{Ratz2006}, the authors used $p = 2$ and $\eta = 30$. In this case,
    \[
\int_{-\infty}^\infty G( U_0)\partial_z  U_0 dz =  30\int_{0}^1  U_0^2(1-U^0)^2 \, d U_0 =    1.
    \]
Let us generalize this result. For any $p \ge 0$, we have
    \[
\int_{-\infty}^\infty G( U_0)\partial_z  U_0 dz = \eta\int_0^1 | U_0|^p|1- U_0|^p\, d U_0 = \eta \frac{\left(\Gamma(1+p)\right)^2}{\Gamma(2+2p)} .    
    \]
A natural choice is the \emph{diffusion normalization},
    \[
\eta =\eta_\star(p) :=  \frac{\Gamma(2+2p)}{\left(\Gamma(1+p)\right)^2}, \quad p \ge 0,
    \]
which guarantees that, for any $p\ge 0$
   \[
 \int_{-\infty}^\infty G( U_0)\partial_z  U_0 dz = 1 .  
    \]
Some common values of the diffusion normalization can be found in Table~\ref{tab:diff-normal}. In conclusion,
	\begin{equation}
W_0 = - \kappa  ,
	\label{chemical-pot-equation-Model-V}
	\end{equation}
which, when combined with the equation
    \[
-v = \partial_{ss}W_0,
    \]
the derivation for which is the same as for Model V and omitted for brevity, leads to the surface diffusion law $v = \partial_{ss} \kappa$.

    \begin{table}[]
    \begin{center}
    \begin{tabular}{c|cccccc}
$p$                  & $\frac{1}{2}$   & 1 & $\frac{3}{2}$         & 2   & 3     & 4 \\ \hline
$\eta_\star(p)$ & $\frac{8}{\pi}$ & 6 & $\frac{128}{3\pi}$ & 30 & 140 & 630
    \end{tabular}
    \end{center}
\caption{Some values of the diffusion normalization $\eta_\star$.}
    \label{tab:diff-normal}
    \end{table}


\begin{thebibliography}{10}

\bibitem{Thompson_ARMR_2012}
Thompson C.~V.. {Solid-state dewetting of thin films}.  {\it Annu. Rev. Mater.
  Res.. }2012;42:399--434.

\bibitem{Baenschetal_JCP_2005}
B\"ansch E., Morin P., Nochetto R.H.. {A finite element method for surface
  diffusion: the parametric case}.  {\it J. Comput. Phys.. }2005;203:321-343.

\bibitem{Hausseretal_IFB_2005}
Hau{\ss}er F., Voigt A.. {A discrete scheme for regularized anisotropic surface
  diffusion: a 6th order geometric evolution equation}.  {\it Interf. Free
  Bound.. }2005;7:353-369.

\bibitem{Baoetal_JCP_2017}
Bao W., Jiang W., Wang Y., Zhao Q.. {A parametric finite element method for
  solid-state dewetting problems with anisotropic surface energies}.  {\it J.
  Comput. Phys.. }2017;330:380-400.

\bibitem{Barrettetal_JCP_2019}
Barrett J.W., Garcke H., N\"urnberg R.. Finite element methods for fourth order
  axisymmetric geometric evolution equations.  {\it J. Comput. Phys..
  }2019;376:733--766.

\bibitem{Wise2005}
Wise S.~M., Lowengrub J.~S., Kim J.~S., Thornton K., Voorhees P.~W., Johnson
  W.~C.. Quantum dot formation on a strain-patterned epitaxial thin film.  {\it
  Appl. Phys. Lett.. }2005;87(13):133102.

\bibitem{Ratz2006}
R{\"a}tz A., Ribalta A., Voigt A.. Surface evolution of elastically stressed
  films under deposition by a diffuse interface model.  {\it J. Comput. Phys..
  }2006;214(1):187--208.

\bibitem{Yeon2006}
Yeon D.-H., Cha P.-R., Grant M.. Phase field model of stress-induced surface
  instabilities: Surface diffusion.  {\it Acta Mater.. }2006;54(6):1623-1630.

\bibitem{Torabi2009}
Torabi S., Lowengrub J., Voigt A., Wise S.. {A new phase-field model for
  strongly anisotropic systems}.  {\it Proc. Royal. Soc. A.
  }2009;465:1337--1359.

\bibitem{Li2009}
Li~B., Lowengrub J., Ratz A., Voigt A.. Geometric evolution laws for thin
  crystalline films: Modeling and numerics.  {\it Commun. Comput. Phys..
  }2009;6(3):433.

\bibitem{Banas2009}
Ba{\v{n}}as L., N{\"u}rnberg R.. Phase field computations for surface diffusion
  and void electromigration in $\mathbb{R}^3$.  {\it Comput. Vis. Sci..
  }2009;12(7):319--327.

\bibitem{Jiang2012}
Jiang W., Bao W., Thompson C.~V, Srolovitz D.~J.. {Phase field approach for
  simulating solid-state dewetting problems}.  {\it Acta Mater..
  }2012;60(15):5578--5592.

\bibitem{Salvalaglio2015a}
Salvalaglio M., Backofen R., Bergamaschini R., Montalenti F., Voigt A..
  {Faceting of equilibrium and metastable nanostructures: A phase-field model
  of surface diffusion tackling realistic shapes}.  {\it Crys. Growth Des..
  }2015;15:2787--2794.

\bibitem{Bergamaschini2016}
Bergamaschini R., Salvalaglio M., Backofen R., Voigt A., Montalenti F..
  {Continuum modelling of semiconductor heteroepitaxy: an applied perspective}.
   {\it Adv. Phys. X. }2016;1(3):331--367.

\bibitem{Naffouti2017}
Naffouti M., Backofen R., Salvalaglio M., et al. Complex dewetting scenarios of
  ultrathin silicon films for large-scale nanoarchitectures.  {\it Sci. Adv..
  }2017;3(11):eaao1472.

\bibitem{Schiedung2017}
Schiedung R., Kamachali R.~D., Steinbach I., Varnik F.. Multi-phase-field model
  for surface and phase-boundary diffusion.  {\it Phys. Rev. E.
  }2017;96:012801.

\bibitem{cahn1996}
Cahn J.W., Elliott C.M., Novick-Cohen A.. The {Cahn-Hilliard} equation with a
  concentration dependent mobility: Motion by minus the {Laplacian} of the mean
  curvature.  {\it Euro. J. Appl. Math.. }1996;7:287-301.

\bibitem{Gugenberger2008}
Gugenberger C., Spatschek R., Kassner K.. Comparison of phase-field models for
  surface diffusion.  {\it Phys. Rev. E. }2008;78:016703.

\bibitem{Dziwniketal_Non_2017}
Dziwnik M., M\"unch A., Wagner A.. {An anisotropic phase-field model for
  solid-state dewetting and its sharp-interface limit}.  {\it Nonlinearity.
  }2017;30:1465-1496.

\bibitem{Voigt2016}
Voigt A.. {Comment on "Degenerate mobilities in phase field models are
  insufficient to capture surface diffusion" [Appl. Phys. Lett. 107, 081603
  (2015)]}.  {\it Appl. Phys. Lett.. }2016;108(3):036101.

\bibitem{Daietal_MMS_2014}
Dai S., Du~Q.. Coarsening mechanism for systems governed by the {Cahn-Hilliard}
  equation with degenerate diffusion mobility.  {\it Multiscale Model. Simul..
  }2014;12:1870-1889.

\bibitem{Leeetal_APL_2015}
Lee A.A., M\"unsch A., S\"uli E.. Degenerate mobilities in phase field models
  are insufficient to capture surface diffusion.  {\it Appl. Phys. Lett..
  }2015;107:081603.

\bibitem{Leeetal_SIAMJAM_2016}
Lee A.A., M\"unsch A., S\"uli E.. Sharp-interface limits of the Cahn-Hilliard
  equation with degenerate mobility.  {\it SIAM J. Appl. Math..
  }2016;76:433-456.

\bibitem{Bhateetal_JAP_2000}
Bhate D.N., Kumar A., Bower A.F.. Diffuse interface model for electromigration
  and stress voiding.  {\it J. Appl. Phys.. }2000;87:1712-1721.

\bibitem{Wiseetal_JCP_2006}
Wise S.~M., Kim J., Lowengrub J.S.. Solving the regularized, strongly
  anisotropic {Cahn-Hilliard} equation by a adaptive nonlinear multigrid
  method.  {\it J. Comput. Phys.. }2006;226:414-446.

\bibitem{dai2012}
Dai S., Du~Q.. Motion of interfaces governed by the {Cahn-Hilliard} equation
  with highly disparate diffusion mobilty.  {\it SIAM J. Appl. Math..
  }2012;72:1818-1841.

\bibitem{dai2016}
Dai S., Du~Q.. Weak solutions for the Cahn-Hilliard equation with degenerate
  mobility.  {\it Arch. Rational Mech. Anal.. }2016;219:1161-1184.

\bibitem{dai2016JCP}
Dai S., Du~Q.. Computational studies of coarsening rates for the
  {Cahn-Hilliard} equation with phase-dependent diffusion mobility.  {\it J.
  Comput. Phys.. }2016;310:85-108.

\bibitem{Karma1998}
Karma A., Rappel W.-J.. Quantitative phase-field modeling of dendritic growth
  in two and three dimensions.  {\it Phys. Rev. E. }1998;57:4323--4349.

\bibitem{Backofen2019}
Backofen R., Wise S.~M., Salvalaglio M., Voigt A.. {Convexity splitting in a
  phase field model for surface diffusion}.  {\it Int. J. Num. Anal. Mod..
  }2019;16(2):192--209.

\bibitem{Albani2016}
Albani M., Bergamaschini R., Montalenti F.. {Dynamics of pit filling in
  heteroepitaxy via phase-field simulations}.  {\it Phys. Rev. B.
  }2016;94(7):075303.

\bibitem{Salvalaglio2015}
Salvalaglio M., Bergamaschini R., Isa F., et al. {Engineered coalescence by
  annealing 3D Ge microstructures into high-quality suspended layers on Si}.
  {\it ACS Appl. Mater. Interfaces. }2015;7(34):19219--19225.

\bibitem{SalvalaglioPRB2016}
Salvalaglio M., Backofen R., Voigt A.. {Thin-film growth dynamics with
  shadowing effects by a phase-field approach}.  {\it Phys. Rev. B.
  }2016;94(23):235432.

\bibitem{Salvalaglio2017a}
Salvalaglio M., Bergamaschini R., Backofen R., Voigt A., Montalenti F., Miglio
  L.. {Phase-field simulations of faceted Ge/Si-crystal arrays, merging into a
  suspended film}.  {\it Appl. Surf. Sci.. }2017;391:33--38.

\bibitem{Salvalaglio2017b}
Salvalaglio M., Backofen R., Voigt A., Montalenti F.. {Morphological evolution
  of pit-patterned Si(001) substrates driven by surface-energy reduction}.
  {\it Nanoscale Res. Lett.. }2017;12(1):554.

\bibitem{Geslin2019}
Geslin P.-A., Buchet M., Wada T., Kato H.. Phase-field investigation of the
  coarsening of porous structures by surface diffusion.  {\it Phys. Rev.
  Materials. }2019;3:083401.

\bibitem{Albani2019}
Albani M., Bergamaschini R., Salvalaglio M., Voigt A., Miglio L., Montalenti
  F.. {Competition between kinetics and thermodynamics during the growth of
  faceted crystal by phase field modeling}.  {\it Phys. Status Solidi (b).
  }2019;256(7):1800518.

\bibitem{Mullins1957}
Mullins W~W. {Theory of thermal grooving}.  {\it J. Appl. Phys.. }1957;28:333.

\bibitem{Mullins1959}
Mullins W.~W.. {Flattening of a nearly plane solid surface due to capillarity}.
   {\it J. Appl. Phys.. }1959;30(1):77--83.

\bibitem{dong2019}
Dong L., Wang C., Zhang H., Zhang Z.. A positivity-preserving, energy stable
  and convergent numerical scheme for the {Cahn-Hilliard} equation with a
  {Flory-Huggins-deGennes} energy.  {\it Commun. Math. Sci..
  }2019;(accepted):1-19.

\bibitem{li2017}
Li~X., Qiao Z., Zhang H.. A second-order convex-splitting Scheme for a
  {Cahn-Hilliard} equation with variable interfacial parameters.  {\it J.
  Comput. Math.. }2017;35:693-710.

\bibitem{Vey2007}
Vey S., Voigt A.. {AMDiS: Adaptive MultiDimensional Simulations}.  {\it Comput.
  Visual Sci.. }2007;10(1):57-67.

\bibitem{Witkowski2015}
Witkowski T., Ling S., Praetorius S., Voigt A.. Software concepts and numerical
  algorithms for a scalable adaptive parallel finite element method.  {\it Adv.
  in Comput. Math.. }2015;41:1145.

\bibitem{HausserJSC2007}
Hau{\ss}er F., Voigt A.. {A discrete scheme for parametric anisotropic surface
  diffusion}.  {\it J. Sci. Comput.. }2007;30:223--235.

\bibitem{BarrettNUMMAT2008}
Barrett J.~W, Garcke H., N{\"{u}}rnberg R.. {Numerical approximation of
  anisotropic geometric evolution equations in the plane}.  {\it Numer. Math..
  }2008;109:1--44.

\bibitem{WangPRB2015}
Wang Y., Jiang W., Bao W., Srolovitz D.~J. {Sharp interface model for
  solid-state dewetting problems with weakly anisotropic surface energies}.
  {\it Phys. Rev. B. }2015;91:045303--045315.

\bibitem{Farina1997}
Farina P.. {The Camunian Rose, Valcamonica Rock Art}.  {\it TRACCE, Online Rock
  Art Bullettin. }1997;7.

\bibitem{SalvalaglioDDCHaniso}
Salvalaglio M., Selch M., Voigt A., Wise S.~M. Doubly degenerate diffuse
  interface models of anisotropic surface diffusion.  {\it Math. Method Appl.
  Sci., in press}. DOI: 10.1002/mma.7118

\end{thebibliography}

\end{document}